%% file: main.tex
\documentclass{article}
\usepackage[utf8]{inputenc}
\input{Packages}

\usepackage{authblk}

\begin{document}

\input{Frontpage}
\newpage
\input{1-Intro}
\input{2-Model} 
\input{3-Tight_approx} 
\input{4-Analysis}
\input{5-Conclusion}
\input{6-Appendix}

\printbibliography

\end{document}

%% file: Packages.tex
\usepackage{graphicx}
\usepackage{wrapfig}
\usepackage{subcaption}
\usepackage{float}
\usepackage[hidelinks]{hyperref}
\usepackage{url}

\usepackage[margin = 2.5 cm]{geometry}

\usepackage[utf8]{inputenc}
\usepackage{gensymb}
\usepackage[english]{babel}

\usepackage{comment}

\usepackage{listings}

\usepackage{amsmath}
\usepackage{mathtools}

\usepackage[backend=biber, style=ieee,doi=false,isbn=false]{biblatex}
\bibliography{References.bib}

%% file: Frontpage.tex
\title{\textbf{Modeling Gas Flow Directions as State Variables: Does it Provide More Flexibility to Power Systems?}}

\author{Junesoo Shin$^{*}$, Yannick Werner$^{*\dagger}$, and Jalal Kazempour$^{*}$ \\
$^{*}$ \small{Department of Wind and Energy Systems, \\ Technical University of Denmark, Kgs. Lyngby, Denmark \\ 
$^{\dagger}$ Department of Industrial Economics and Technology Management, \\
Norwegian University of Science and Technology,
Trondheim, Norway\\
s184386@student.dtu.dk, $\{$yanwe, seykaz$\}$@dtu.dk}
}


\maketitle
\vspace{20pt}
\begin{abstract}
As a common practice, the direction of natural gas flow in every pipeline is determined ex-ante for simplification purposes, and treated as a given parameter within the scheduling problem. However, in  integrated gas and electric power networks with a large share of intermittent renewable power supply, it is no longer straightforward to optimally predetermine the gas flow directions. A wrong predetermination of gas flow directions may result in feasible but not necessarily optimal schedules. We propose a mixed-integer linear optimization model to determine the optimal gas flow directions while scheduling the system. This unlocks additional flexibility to power systems, provided that  a  tight coordination between power and gas systems exists. The increased flexibility, although it comes at the cost of increased computational complexity, is quantified by comparing the total operational cost of the entire system with bidirectional gas flows as opposed to unidirectional gas flows. We numerically show that modeling gas flow directions as state variables may bring added value not only in the meshed but also in the radial gas networks. \vspace{20pt}\\
\textbf{Keywords:} Integrated energy systems, optimal gas flow directions, flexibility, linepack, mixed-integer linear program.
\end{abstract}

%% file: 1-Intro.tex
\section{Introduction}
\label{Section 1}
\subsection{Motivation and aim}
The increasing share of intermittent renewable power supply increases the need for operational flexibility in electric power systems. Natural gas-fired power plants are generally able to provide flexibility by quickly adjusting their power output in order to cope with fluctuations in the renewable power supply \cite{gil}. This, however, leads to fluctuations in the natural gas demand and might even trigger a change of the flow direction in certain gas pipelines. If this change contradicts the predetermined flow direction, the operational space of natural gas-fired power plants will be shrunk. In this case, more expensive flexibility options have to be utilized, leading to an increase of total operational cost of the system.

Pursuing modeling simplicity, it is currently a common practice to assume that flow directions in natural gas pipelines can be determined ex-ante based on historical observations and expert judgement, no matter what the schedule of gas-fired power plants is \cite{Massol2018MarketHubs}. This predetermined flow direction is then treated as a fixed parameter within the scheduling problem. 
This assumption is  more frequently taken in gas networks with a radial topology, as it is the case in Denmark, compared to meshed gas networks, which is the case in Belgium \cite{EuropeanPlatform}. Nonetheless, it has been argued that flow directions even in meshed networks can be predetermined in a straightforward manner, e.g., by using a simplified model  neglecting the gas flow dynamics in pipelines \cite{Munoz2003NaturalStudies}. 

This practice is, however, being challenged, as the uncertainty induced by renewable power supply propagates from the power system to the natural gas network via coupling components at the interface of the two systems, in particular natural gas-fired power plants. As a consequence, natural gas demands and eventually optimal flow quantities in pipelines and even flow directions become uncertain \cite{Keyaerts2014ImpactEurope}. One may hypothesize this uncertainty will be even more magnified in systems with significant stochastic renewable gas injection. 

This paper argues that in order to unlock the full flexibility potential that the natural gas network can provide to the power system, not only the schedules and flow quantities but also the gas flow directions should be optimally determined. This calls for a change in the current practice that predetermines the gas flow directions. We find out that even in radial gas networks it might no longer be straightforward to optimally predetermine the flow directions, in particular in cases that natural gas suppliers are located far away from each other throughout the network. 

Hereafter, by \textit{unidirectional} we refer to those models with predetermined gas flow directions, while by \textit{bidirectional} we refer to models where the gas flow direction in every pipeline is optimally determined along with schedules and flow quantities. 

\subsection{State of the art}

\begin{table}[h]
\centering
\caption{State of the art}
\resizebox{0.7\textwidth}{!}{%
\begin{tabular}{lllll}
Reference & Convexification method & Model & \begin{tabular}[c]{@{}l@{}}Bidirectional \\ gas flow  modeling \end{tabular} & \begin{tabular}[c]{@{}l@{}}Linepack  \\ modeling\end{tabular} \\ 
\hline \vspace{5pt}
\cite{D.Manshadi2019CoordinatedApproach}& Convex relaxation & SDP & Yes & No \\
\vspace{5pt}

\cite{He2018DecentralizedReformulation} & Convex relaxation & MISOCP & Yes & No \\
\vspace{5pt} 

\cite{Chen2019UnitModel} & Convex relaxation & MISOCP & No & Yes  \\
\vspace{5pt} 

\cite{Schwele2019a} & Convex relaxation & MISOCP & Yes & Yes \\
\vspace{5pt}

\cite{Sirvent2017LinearizedSystems}, \cite{Munoz2003NaturalStudies}& Piecewise linear approximation & MILP & Yes & No \\
\vspace{5pt}

\cite{Correa-Posada2015IntegratedOperation} & Piecewise linear approximation & MILP & Yes & Yes \\
\vspace{5pt}

\cite{Ordoudis2019} & Outer linear approximation & MILP & Yes & Yes \\
\vspace{5pt}

This paper & \begin{tabular}[c]{@{}l@{}}Outer linear approximation \\ (tighter constraints)\end{tabular} & MILP & Yes & Yes
\end{tabular}%
}
\label{Table_lit_review}
\vspace{-0.5cm}
\end{table}

The current literature proposes a broad range of operational models for integrated power and natural gas systems, mainly built upon convexified gas flow equations, while taking into account the bidirectionality of gas flows and linepack (stored gas in the pipelines). Table \ref{Table_lit_review} provides an overview. The interested reader is referred to \cite{Raheli2021OptimalMethods} for a more comprehensive  survey.
Reference \cite{D.Manshadi2019CoordinatedApproach} proposes a semidefinite programming (SDP) approach, where a tractable and tight convex relaxation is explored using many small matrices with lifting variables. However, for dense gas networks, it becomes challenging to solve the resulting model due to large matrices.
In addition, \cite{He2018DecentralizedReformulation}, \cite{Chen2019UnitModel}, and \cite{Schwele2019a} use a second-order cone (SOC) relaxation to convexify the gas flow equations, while proposing solutions to tighten the relaxation. 
Reference \cite{He2018DecentralizedReformulation} proposes an iterative algorithm based on sequential cone programming to obtain a reliable solution to the original problem. 
Similarly, \cite{Chen2019UnitModel} develops an iterative algorithm that is tightening the bounds of the McCormick envelops used to convexify the bilinear terms that are necessary to convexify the gas flow equations when accounting for linepack. 
Reference \cite{Schwele2019a} investigates both linepack modeling and bidirectionality, while proposing a mixed-integer linear program (MILP) as well as a mixed-integer second-order cone program (MISOCP). To convexify bilinear terms, a McCormick relaxation is used. 
References \cite{Munoz2003NaturalStudies}, \cite{Sirvent2017LinearizedSystems}, and \cite{Correa-Posada2015IntegratedOperation} propose piecewise linear approximation methods to convexify the gas flow equations. In detail,  \cite{Munoz2003NaturalStudies} proposes a two-stage model, where the flow direction of the passive pipelines is determined and then a non-linear program is solved.
In addition, \cite{Sirvent2017LinearizedSystems} develops an extended incremental method including a finitely bounded variable.
Reference \cite{Correa-Posada2015IntegratedOperation} uses an incremental formulation of the piecewise linear approximation, as originally suggested by \cite{Correa-Posada2014GasModels}.
Finally, \cite{Ordoudis2019} exploits an outer linear approximation of the gas flow equations, while exploring sequential and stochastic coupling of power and gas systems. Due to its simplicity and computational advantages compared to other approaches presented here, we use an outer linear approximation based on a Taylor-series expansion. For that, we use a set of predefined fixed pressure points to approximate non-convex gas flow dynamics.

\subsection{Contributions}
We quantify the increased flexibility potential that the natural gas system provides to the power system by including the gas flow directions as state variables in the scheduling problem. We consider a co-optimization problem for the integrated power and gas system, accounting for linepack in the natural gas pipelines. After some reformulations to achieve convexity, we compare two cases: The first case includes only nodal pressures and pipeline flows as state variables and assumes predetermined gas flows (unidirectional). In the second case, the flow directions are included as state variables in the optimization model (bidirectional). This brings an extra degree of freedom to the system operator but at the potential expense of increased computational complexity due to the inclusion of binary variables for representing the flow directions.
Using a stylized case study, we provide an in-depth analysis of natural gas flow quantities, directions, and linepack utiliziation in both cases considering a meshed as well as a radial gas network. We show that including gas flow directions as state variables may provide substantial flexibility to the power system and foster the utilization of linepack, leading to a decrease in total operational cost of the integrated system. Methodologically, we contribute to the literature by improving the currently prevalent linear approximation of the non-linear and non-convex Weymouth equation that governs the gas flow dynamics by further tightening the feasible region of the approximation. This ensures that the optimal flow directions obtained with the convexified model remain the same when recalculating them based on the optimal pressures using the original Weymouth equation.

\vspace{-0.1cm}
\subsection {Paper organization}
The remainder of this paper is organized as follows. Section~\ref{Section 2} introduces optimization models for the operation of integrated power and gas system with unidirectional and bidirectional flows. Section~\ref{Section 3} proposes a tightened version of the outer linear approximation of the non-convex gas flow equations for both unidirectional and bidirectional cases. Section~\ref{Section 4} presents a detailed analysis of both models applied to a stylized case study, considering both  meshed and radial gas network topologies. Section~\ref{sec: Section 5} concludes. Finally, Appendices A-C provide further modeling details. 

%% file: 2-Model.tex
\vspace{-0.1cm}
\section{Formulation}
\label{Section 2}
\subsection{General model}
Aiming at harnessing further flexibility from the gas network by modeling flow directions as state variables, we consider an integrated power and natural gas system given perfect and instantaneous information exchange. Although such a fully coordinated system does not represent the current status of real-world system operation, it provides an ideal benchmark to evaluate the maximum flexibility potential in the coordinated system when accounting for bidirectionality of gas flows. We consider a deterministic model with a single-point forecast and inelastic demand profiles while optimally determining dispatch schedules of power generators and gas suppliers in the day-ahead time stage. A compact form of the co-optimization model is given below:
\begin{subequations}
\begin{align}
    \underset{\textbf{x},\textbf{y}}{\operatorname{Minimize}} \quad &Cost^{\text{power}}(\mathbf{x}) + Cost^{\text{gas}}(\mathbf{y}) \label{compact1}\\
    & \mathbf{g}^{\text{power}}(\mathbf{x}) \leq 0 \label{compact2}\\
    & \mathbf{h}^{\text{power}}(\mathbf{x}) = 0 \label{compact3}\\
    & \mathbf{g}^{\text{gas}}(\mathbf{y}) \leq 0 \label{compact4}\\
    & \mathbf{h}^{\text{gas}}(\mathbf{y}) = 0 \label{compact5} \\
    & \mathbf{e}^{\text{power}}(\mathbf{x}) + \mathbf{f}^{\text{gas}}(\mathbf{y}) = 0. \label{compact6}
\end{align}
\end{subequations}

The objective function \eqref{compact1} minimizes the total operational cost of the integrated power and gas system. The capacity, flow, and nodal balance constraints for power and gas systems are represented as \eqref{compact2}-\eqref{compact6}, where \eqref{compact6} links the power and gas systems together. The power and gas system variables are represented by vectors $\mathbf{x}$ and $\mathbf{y}$, respectively. Throughout this paper, upper-case letters denote parameters whereas variables are represented by lower-case letters.

Since the focus of this work lies on the difference between unidirectional and bidirectional gas flow models, we present the objective function \eqref{compact1} as well as the power system  constraints \eqref{compact2}-\eqref{compact3}, which are shared in both models, in Appendix A and Appendix B, respectively.  Section~\ref{Common gas system constraints} describes those gas network constraints which are common in both unidirectional and bidirectional models. Finally, we  provide a detailed description of the gas flow constraints \eqref{compact4}-\eqref{compact5} as well as linking constraints \eqref{compact6} accounting for unidirectional and bidirectional flows in Section~\ref{Unidirectional gas flow contraints} and Section~\ref{Bidirectional gas flow contraints}, respectively. 

\subsection{Common gas system constraints}
\label{Common gas system constraints}
Let $m, u \in \mathcal{N}$, $k \in \mathcal {K}$, and $t \in \mathcal{T}$ denote the set of gas network nodes, gas suppliers, and time periods, respectively. Each gas pipeline connects two nodes of the network, such that $(m,u) \in \mathcal{Z}$ denotes the set of all pipelines in the network.
In order to properly model technical limits of the gas network, the following set of constraints is enforced:
\begin{subequations}\label{common_gas_constraints}
\begin{align}
    &0 \leq g_{k,t} \leq G^{\text{max}}_k , \hspace{2mm}\forall k,t \label{g1}\\
    &\mathit{PR}^{\text{min}}_m \leq pr_{m,t} \leq \mathit{PR}^{\text{max}}_m, \hspace{2mm}\forall m,t \label{g2}\\
    &pr_{u,t} \leq \Gamma_{m,u} pr_{m,t}, \hspace{2mm}\forall (m,u) \in \mathcal{Z},t. \label{g3}
\end{align}
\end{subequations}

Constraints \eqref{g1} limit the schedule $g_{k,t}$ of gas supply unit $k$ in time period $t$ to its maximum capacity $G^{\text{max}}_k$. Constraints \eqref{g2} enforce the nodal pressures $pr_{m,t}$ at node $m$ in time period $t$ to lie within the technical limits $\mathit{PR}^{\text{min}}_m$ and $\mathit{PR}^{\text{max}}_m$. We use a simplified representation of compressors \eqref{g3} in pipeline $(m,u) \in \mathcal{Z}$, where we assume a constant compression ratio $\Gamma_{m,u} > 1$ for pipelines that host compressors and $\Gamma_{m,u} = 1$ otherwise. Furthermore, any form of energy demand from the compressor is neglected.

\subsection{Unidirectional gas flow constraints}
\label{Unidirectional gas flow contraints}
The steady-state gas flow $q_{m, u, t}$ in pipeline $(m,u)$ and time period $t$ is determined by the Weymouth equation as
\begin{subequations}
\begin{flalign}
    &q_{m,u,t}^2 = K^{2}_{m,u} (pr^{2}_{m,t} - pr^{2}_{u,t}), \hspace{2mm}\forall (m,u) \in \mathcal{Z},t,\label{eq:Weymouth}
\end{flalign}
\end{subequations}
where $K_{m,u}$ denotes the natural gas flow constant of pipeline $(m,u)$, capturing the technical characteristics of the pipeline. In the unidirectional model, we assume that the flow direction is predetermined and fixed to flow from node $m$ to node $u$ in all time periods, which is enforced by 
\begin{subequations}\label{LP_qs}
\begin{gather}
    q_{m,u,t} \geq 0 \hspace{2mm}\forall(m, u) \in \mathcal{Z}, t, \label{fixed direction}
\end{gather}
In addition, we enforce 
\begin{gather}
    q_{m, u, t}=\frac{q_{m, u, t}^{\text {in }}+q_{m, u, t}^{\text {out}}}{2}, \hspace{2mm} \forall(m, u) \in \mathcal{Z}, t. \label{LP:q1}
\end{gather}
\end{subequations}
which relates the gas flow within a pipeline to the average of its gas inflow $q_{m, u, t}^{\text {in}}$ and outflow $q_{m, u, t}^{\text {out}}$.

Due to the slow transients of natural gas flows in pipelines, a time delay between gas inflow and outflow exists, which behaves like an intertemporal storage. This storage can be described by the following set of constraints:
\begin{subequations}\label{LP_linepack_constraints}
\begin{align}
    &h_{m, u, t}=S_{m, u} \frac{p r_{m, t}+p r_{u, t}}{2}, \hspace{2mm}\forall(m, u) \in \mathcal{Z}, t \label{LP:h1} \\
    &h_{m, u, t}=h_{m, u,(t-1)}+q_{m, u, t}^{\text {in }}-q_{m, u, t}^{\text {out }}, \forall(m, u) \in \mathcal{Z}, t \label{LP:h2} \\
    &H_{m, u}^{0} \leq h_{m, u, t}, \hspace{2mm}\forall(m, u) \in \mathcal{Z}, t=|\mathcal{T}|. \label{LP:h4}
\end{align}
\end{subequations}
Constraints \eqref{LP:h1} relate the linepack mass $h_{m, u, t}$ in the pipeline $(m,u)$ and time period $t$ to the average pressure difference of the connected nodes and the pipeline's specific linepack constant $S_{m, u}$, which, again, captures the technical characteristics of the pipeline. The intertemporal energy balance is defined by \eqref{LP:h2}, where the initial linepack mass $h_{m, u, 0}$ is fixed to a given level $H_{m, u}^{0}$, by setting $h_{m, u, 0} = H_{m, u}^{0}$. Constraints \eqref{LP:h4} ensure that a depletion of natural gas in the pipelines at the end of the planning horizon, i.e.,  $t=|\mathcal{T}|$, is avoided by lower bounding the linepack mass to the initial level.

As linking constraints in the form of \eqref{compact6}, the nodal  balance for every gas node $m$ is enforced by
\begin{gather}
    \sum_{k \in \mathcal{A}_{m}^{\rm{K}}} g_{k, t}-\sum_{i \in \mathcal{A}_{m}^{\rm{G}}} \eta_{i} p_{i, t} -\sum_{u:(m, u) \in \mathcal{Z}}\left(q_{m,u,t}^{\text{in}}-q_{m,u,t}^{\text{out}}\right)
    \hspace{0.5cm}=\sum_{d \in \mathcal{A}^{\rm{DG}}_m} D_{d, t}^{\mathrm{G}}, \hspace{2mm} \forall m, t. \label{LP_gas_bal}
\end{gather}
where $\mathcal{A}_{m}^{\rm{K}}$, $\mathcal{A}_{m}^{\rm{G}}$, and $\mathcal{A}^{\rm{DG}}_m$ describe the sets of gas supply units, gas-fired power plants (GFPPs), and gas loads in node $m$, respectively.
We consider the gas supply schedule $g_{k, t}$ and power production schedule $p_{i,t}$ by GFFP $i$ to be variable, while we assume that the gas demand $D_{d, t}^{\mathrm{G}}$ of every gas demand $d$ is fixed. The power conversion factor of GFPPs, $\eta_{i}$, links the power and gas system operation by relating the gas demand of GFFP $i$ to its power production. 

The unidirectional gas flow constraints constitute a linear problem without binary variables. In the following section, we  show how this model can be extended to account for bidirectional flows, resulting in a mixed-integer linear program.

\subsection{Bidirectional gas flow constraints}
\label{Bidirectional gas flow contraints}
In order to account for bidirectional flows, the unidirectional Weymouth equation \eqref{eq:Weymouth} is replaced by
\begin{subequations}
\begin{flalign}
    &q_{m,u,t} |q_{m,u,t}| = K^{2}_{m,u} (pr^{2}_{m,t} - pr^{2}_{u,t}), \hspace{2mm}\forall (m,u) \in \mathcal{Z},t,\label{eq:MILP_Weymouth}
\end{flalign}
\end{subequations}
which allows the gas flow $q_{m,u,t}$ to be either positive, or zero, or negative, depending on the sign of the pressure difference $(pr^{2}_{m} - pr^{2}_{u})$ at the adjacent nodes. 

Similarly, constraints \eqref{LP_qs} are replaced by the following set of constraints:
\begin{subequations}\label{MILP_qs}
\begin{align}
    &q_{m,u,t} = q^+_{m,u,t} - q^-_{m,u,t}, \hspace{2mm}\forall (m,u) \in \mathcal{Z},t \label{MILP:q1}\\
    &0 \leq q_{m, u, t}^{+} \leq M y_{m, u, t}, \hspace{2mm}\forall(m, u) \in \mathcal{Z}, t \label{MILP:q2}\\
    &0 \leq q_{m, u, t}^{-} \leq M\left(1-y_{m, u, t}\right), \hspace{2mm}\forall(m, u) \in \mathcal{Z}, t \label{MILP:q3}\\
    &q_{m, u, t}^{+}=\frac{q_{m, u, t}^{\text {in }}+q_{m, u, t}^{\text {out }}}{2}, \hspace{2mm}\forall(m, u) \in \mathcal{Z}, t \label{MILP:q4}\\
    &q_{m, u, t}^{-}=\frac{q_{u, m, t}^{\text {in }}+q_{u, m, t}^{\text {out }}}{2}, \hspace{2mm}\forall(m, u) \in \mathcal{Z}, t \label{MILP:q5} \\
    &q_{m,u,t} \geq 0, \hspace{2mm}\forall \{(m,u) \in \mathcal{Z}|\Gamma_{m,u}\neq1\},t \label{MILP:q6}\\
    &y_{m, u, t} \in \{0,1\}, \forall(m, u) \in \mathcal{Z}, t \label{MILP:q7}. 
\end{align}
\end{subequations}
The gas flow $q_{m,u,t}$ is described in \eqref{MILP:q1} by two non-negative variables $q_{m, u, t}^{+}$ and $q_{m, u, t}^{-}$ that represent flow directions from node $m$ to $u$ and node $u$ to $m$, respectively. Using binary variable $y_{m, u, t}$,  \eqref{MILP:q2} and \eqref{MILP:q3} enforce that the gas can only flow in one direction in each pipeline and time period using the big-M method. The value of constant $M$ is chosen to be sufficiently high to not additionally constrain the amount of natural gas flows, but  not too high to avoid causing numerical problems. Constraints \eqref{MILP:q4} and \eqref{MILP:q5} relate the flow in a pipeline to the average of in- and outflow.
In pipelines that host compressors, the flow direction is not governed by the Weymouth equation, since the gas flows from a node with lower pressure to a node with higher pressure. Hence, \eqref{MILP:q6} ensures that the flow in these pipelines is unidirectional. Constraints \eqref{MILP:q7} introduce binary variables. 

Accounting for bidirectional flows in pipelines, the linepack constraints \eqref{LP:h2} have to be adjusted as 
\begin{gather}
    h_{m, u, t}=h_{m, u,(t-1)}+q_{m, u, t}^{\text {in }}-q_{m, u, t}^{\text {out }}+q_{u, m, t}^{\text{in}}-q_{u, m, t}^{\text{out}},
    \forall(m, u) \in \mathcal{Z}, t. \label{MILP:h2}
\end{gather}

Similarly, the nodal gas balance constraints \eqref{LP_gas_bal} are rewritten as
\begin{gather}
    \sum_{k \in \mathcal{A}_{m}^{\rm{K}}} g_{k, t}-\sum_{i \in \mathcal{A}_{m}^{\rm{G}}} \eta_{i} p_{i, t} -\sum_{u:(m, u) \in \mathcal{Z}} \Big(q_{m,u,t}^{\text{in}}-q_{m,u,t}^{\text{out}}  
    \hspace{0.5cm}+q_{u,m,t}^{\text{in}}-q_{u,m,t}^{\text{out}}\Big) =\sum_{d \in \mathcal{A}^{\rm{DG}}_m} D_{d, t}^{\mathrm{G}}, \hspace{2mm} \forall m, t. \label{MILP_gas_bal}
\end{gather}

Note that both unidirectional and bidirectional gas flow models are non-convex due to the Weymouth equation \eqref{eq:Weymouth} and \eqref{eq:MILP_Weymouth}. Hence, we  adopt a convex linear approximation of the Weymouth equation in the next section. Furthermore, we improve the approximation that is currently prevalent in the literature \cite{Ordoudis2019,Tomasgard2007OptimizationChain} by showing how it can be further tightened to account for some of the underlying gas flow physics. 

%% file: 3-Tight_approx.tex
\section{Tighter approximation}
\label{Section 3}
\begin{figure}[H]
\centering
\includegraphics[width=9.5cm,height=6cm]{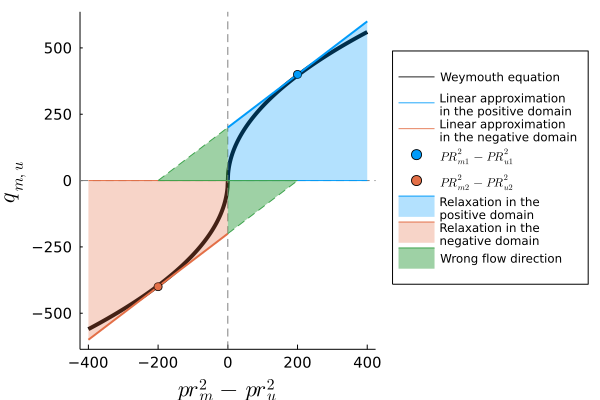}
\caption{Outer linear approximation of the Weymouth equation.}
\label{fig:approx}
\vspace{-0.2cm}
\end{figure}
In Appendix C, Weymouth equation is convexified using a Taylor series approximation and relaxation. This outer approximation method violates some of the underlying physics of gas flows. The main issues are two-fold:
\begin{itemize}
    \item The linear approximations extend to the green areas by crossing the y-axis as shown in Fig. \ref{fig:approx}. The points in the green area treat gas flows from nodes with lower pressure to higher pressure as feasible, since the sign of $pr^{2}_{m} - pr^{2}_{u}$ is flipped.
    \item The relaxation includes points on the y-axis and x-axis. Apart from the origin, points on the y-axis have gas flows when there is no pressure difference in the adjacent nodes. Similarly, points on the x-axis have no gas flow when there is pressure difference in the adjacent nodes.
\end{itemize}
Since determining the optimal gas flow directions is critical to our model, we resolve the above approximation issues by introducing the following additional  constraints that further tighten the feasible space.
\begin{subequations}\label{LP_additional_constraints}
\begin{align}
    &pr_{m,t} \geq pr_{u,t},\hspace{2mm}\forall (m,u) \in \mathcal{Z},t \label{LP_issue1} \\
    &0 \leq q_{m, u, t} \leq M(pr_{m,t}-pr_{u,t}), \hspace{2mm}\forall (m,u) \in \mathcal{Z},t. \label{LP_issue2}
\end{align}
\end{subequations}
Constraints \eqref{LP_issue1} enforce gas flow from nodes with higher pressure to lower pressure in the unidirectional model. To prevent gas flow in the pipelines with zero pressure difference between the adjacent nodes, \eqref{LP_issue2} is introduced. However, the extent to which one value can be considered larger than the other is arbitrary, allowing the gas flow with very small pressure difference to be feasible. The same issue exists for very small gas flow while having a large difference in the pressures in the adjacent nodes. These issues are inherent to any convex relaxation technique and are therefore not further addressed in this paper. 

Constraints \eqref{LP_issue1} are adjusted to account for bidirectional flows in pipelines as 
\begin{subequations}
\begin{align}
    &(pr_{m,t} - pr_{u,t}) y_{m,u,t} \geq 0, \hspace{2mm}\forall (m,u) \in \mathcal{Z},t \label{MILP_issue1} \\
    &(pr_{u,t} - pr_{m,t}) (1-y_{m,u,t}) \geq 0, \hspace{2mm}\forall (m,u) \in \mathcal{Z},t. \label{MILP_issue2}
\end{align}
\end{subequations}
However, \eqref{MILP_issue1} and \eqref{MILP_issue2} are non-linear due to the product of  continuous variables $pr_{m,t}$ and binary variables $y_{m,u,t}$. We use a binary expansion method, that exactly reformulates the above non-linear expressions as linear expressions, as given in \eqref{MILP_additional_constraints}. For the implementation of the binary expansion method, two non-negative auxiliary variables $\phi_{m,t}$ and $\phi_{u,t}$ are introduced, where $\phi_{m,t}=pr_{m,t}y_{m,u,t}$ and $\phi_{u,t}=pr_{u,t}y_{m,u,t}$.
\begin{subequations}\label{MILP_additional_constraints}
\begin{align}
    &\phi_{m,t}  - \phi_{u,t} \geq 0, \hspace{2mm}\forall (m,u) \in \mathcal{Z},t \label{MILP_issue3}\\
    &pr_{u,t} - pr_{m,t} - \phi_{u,t} +\phi_{m,t} \geq 0, \hspace{2mm}\forall (m,u) \in \mathcal{Z},t \label{MILP_issue4} \\
    &-y_{m,u,t} M \leq \phi_{m,t} \leq y_{m,u,t} M, \hspace{2mm}\forall (m,u) \in \mathcal{Z},t \label{MILP_issue5}\\
    &-(1-y_{m,u,t}) M \leq \phi_{m,t}-pr_{m,t} \leq (1-y_{m,u,t})M, \notag \\
    &\hspace{0.5cm} \hspace{2mm}\forall (m,u) \in \mathcal{Z},t \label{MILP_issue6}\\
    &-y_{m,u,t} M \leq \phi_{u,t} \leq y_{m,u,t} M, \hspace{2mm}\forall (m,u) \in \mathcal{Z},t \label{MILP_issue7}\\
    &-(1-y_{m,u,t}) M \leq \phi_{u,t}-pr_{u,t} \leq (1-y_{m,u,t})M, \notag \\
    &\hspace{0.5cm} \hspace{2mm}\forall (m,u) \in \mathcal{Z},t. \label{MILP_issue8}
\end{align}
\end{subequations}
Constraints \eqref{LP_issue2} are adjusted as \eqref{MILP_issue9} and \eqref{MILP_issue10} to account for bidirectional gas flow. 
\begin{subequations}
\begin{align}
    &0 \leq q_{m, u, t}^{+} \leq M(pr_{m,t}-pr_{u,t})y_{m,u,t}, \hspace{2mm} \forall (m,u) \in \mathcal{Z},t \label{MILP_issue9} \\
    &0 \leq q_{m, u, t}^{-} \leq M(pr_{u,t}-pr_{m,t})(1-y_{m,u,t}), \hspace{2mm} \forall (m,u) \in \mathcal{Z},t. \label{MILP_issue10}
\end{align}
\end{subequations}
Dealing with the same non-linearity issue, the auxiliary variables $\phi_{m,t}$ and $\phi_{u,t}$ are introduced to linearize \eqref{MILP_issue9} and \eqref{MILP_issue10} to \eqref{MILP_issue11} and \eqref{MILP_issue12}.
\begin{subequations}\label{MILP_additional_constraints2}
\begin{align}
    &0 \leq q_{m, u, t}^{+} \leq M(\phi_{m,t}-\phi_{u,t}), \hspace{2mm} \forall (m,u) \in \mathcal{Z},t \label{MILP_issue11} \\
    &0 \leq q_{m, u, t}^{-} \leq M(pr_{u,t}-pr_{m,t}-\phi_{u,t}+\phi_{m,t}), \hspace{2mm} \forall (m,u) \in \mathcal{Z},t. \label{MILP_issue12}
\end{align}
\end{subequations}
Table \ref{Table:LPvsMILP} summarizes the structure of optimization problems for both approximated unidirectional and bidirectional gas flow models.

\begin{table}[H]
\centering
\caption{Summary of approximated gas flow models.}
\resizebox{0.7\textwidth}{!}{%
\begin{tabular}{l|cc}
 & Optimization problem & The set of variables \\ \hline
\begin{tabular}[c]{@{}l@{}}Unidirectional\\ (Linear program)\end{tabular} 
& 
\begin{tabular}[c]{@{}l@{}}
Objective function: \eqref{objective} \\
Power system constraints: \eqref{power constraints} \\
Gas system constraints: \\
\eqref{common_gas_constraints}, \eqref{LP:q1}, \eqref{LP_linepack_constraints}, \eqref{LP_gas_bal}, \eqref{LP_additional_constraints}, \eqref{eq:LP_weymouth}
\end{tabular} 
&
\begin{tabular}[c]{@{}l@{}} 
$\Theta_{\mathrm{Unidirectional}} =$ \\ 
$\{p_{i, t}, w_{j, t}, g_{k, t}, \theta_{n, t}, f_{n,r,t},$ \\
$h_{m, u, t},$, $pr_{m, t}, q_{m, u, t}, q_{m, u, t}^{\mathrm{in}}, q_{m, u,t}^{\mathrm{out}}\}$
\end{tabular}\\ \hline
\begin{tabular}[c]{@{}l@{}}Bidirectional\\ (Mixed-integer \\ linear program)\end{tabular} 
& 
\begin{tabular}[c]{@{}l@{}} 
Objective function: \eqref{objective} \\
Power system constraints: \eqref{power constraints} \\
Gas system constraints: \\
\eqref{common_gas_constraints}, 
\eqref{LP:h1}, \eqref{LP:h4}, \eqref{MILP_qs}, 
\eqref{MILP:h2}, 
\eqref{MILP_gas_bal}, \\ \eqref{MILP_additional_constraints},  \eqref{MILP_additional_constraints2}, \eqref{MILP_flow_approx}
\end{tabular} 
& 
\begin{tabular}[c]{@{}l@{}}
$\Theta_{\mathrm{Bidirectional}} =$ \\
$\{ \Theta_{\mathrm{Unidirectional}}, q_{m, u, t}^{+},q_{m, u, t}^{-},$\\
$q_{u, m, t}^{\text {in }}q_{u, m, t}^{\text {out }}, y_{m, u, t}, \phi_{m,t},\phi_{u,t}\}$ \end{tabular}
\label{Table:LPvsMILP}
\end{tabular}}
\end{table}

%% file: 4-Analysis.tex
\section{Case study}
\label{Section 4}
The case study is a stylized $24$-node IEEE reliability test power system, connected to a $12$-node natural gas system that is composed of $12$ power generators ($5$ non-GFPPs and $7$ GFPPs), $2$ wind farms, $3$ gas supply units, $17$ electricity loads, and $4$ gas loads. The online appendix for the original case study with network topology is available at \cite{Ordoudis2019}. For all cases, the total installed wind power capacity was chosen to equal $40$\% of the hourly average electricity demand.
The models are solved using Lenovo SD$530$ with two Intel Xeon Gold $6226$R processors ($2.90$ GHz, $16$ core per CPU), and $378$ GB RAM in Julia using JuMP and Gurobi solver package $9.1.0$. 
Due to computational complexities when solving the mixed-inter linear program that describes the bidirectional gas flow model, the time horizon was split into two parts: hours $1$ to $12$ and hours $13$ to $24$. To avoid a depletion of natural gas in the pipelines, the linepack mass at hour $12$ is used as the initial line pack mass of hour $13$.
\begin{figure}[H]
\centering
\includegraphics[width=16.5cm]{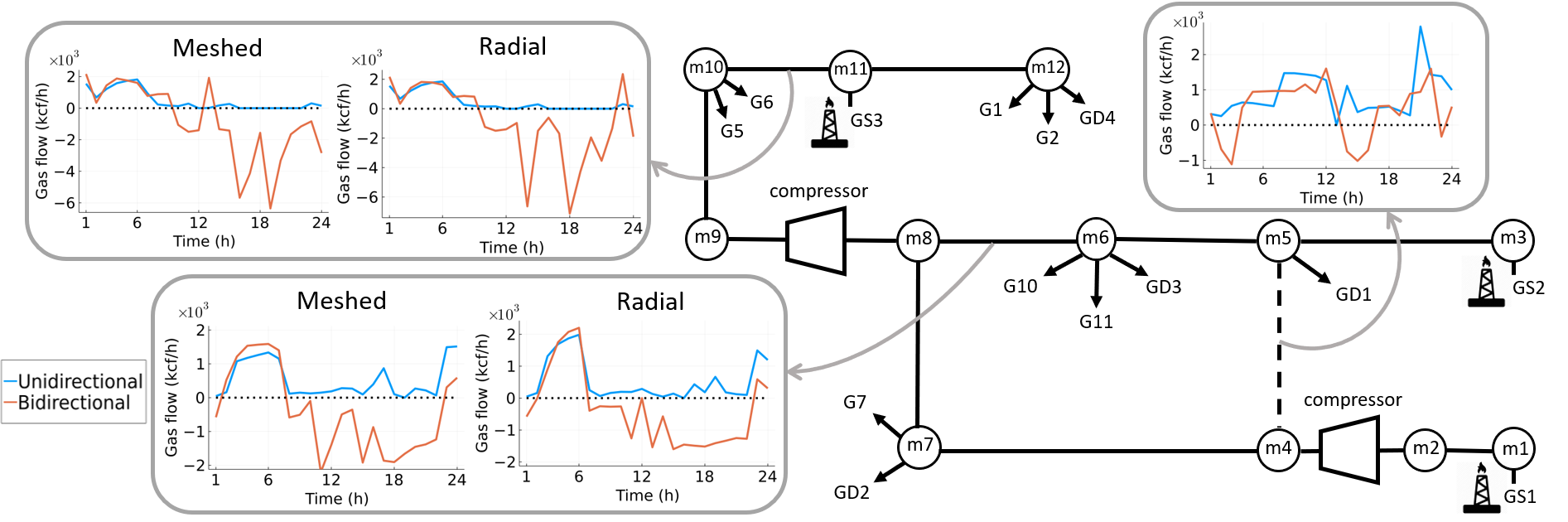}%
\caption{Gas flows for unidirectional and bidirectional models in the pipelines ($m4-m5$), ($m6-m8$), and ($m10-m11$). If the pipeline ($m4-m5$) is present, it represents the meshed gas network, whereas if the pipeline ($m4-m5$) is removed, it represents the radial gas network.}
\label{fig:LPvsMILP}
\end{figure}
To study the benefit of modeling bidirectional gas flows as opposed to unidirectional gas flows in both radial and meshed systems, the gas system network in the case study is slightly adjusted. Fig. \ref{fig:LPvsMILP} shows the gas network diagram, including nodes $m1$ to $m12$, gas supply units GS$1$ to GS$3$, gas loads GD$1$ to GD$3$, and GFPPs G$1$, G$2$, G$5$, etc. The dotted pipeline between nodes $4$ and $5$ in the meshed network is removed in the radial network topology.
Note that the gas suppliers are located at different ends of the network that are topologically far away from each other.

We observed that the unidirectional gas flow model can be solved within seconds for both meshed and radial systems, whereas the bidirectional gas flow model was solved in $2.6$ hours for the meshed system, and in less than two minutes for the radial system. All input data and codes used in this paper are publicly available in the online companion \cite{Shin2022ModelingSystems}. 

\subsection{Model comparison}
For the bidirectional gas flow model, we observe that a change in the flow directions happens most frequently in the pipelines ($m4-m5$), ($m6-m8$), and ($m10-m11$), which are connected to nodes with gas loads and GFPPs. The gas flow in these pipelines is shown in the attached plots in Fig. \ref{fig:LPvsMILP}. It can be further noticed, that the change of flow direction happens not only in the meshed but also in the radial system. Moreover, the changes in flow directions mostly occur in the later half of the day when the demands for electricity and gas are higher. This is accompanied by a higher amount of gas that is transported through the pipelines, indicating an increase of the utilization of comparatively cheap GFPPs. In fact, in the meshed system, the total power production by GFFPs as a share of total electricity demand increases from $13.0$\% in the unidirectional case to $16.7$\% in the bidirectional case. Consequently, allowing for bidirectional flows in the gas systems increases the flexibility in the power system. This leads to an operational cost saving of $2.1$\% in the meshed system and $1.2$\% in the radial system.
The next section describes the difference in gas flows and power production schedule of GFPPs between the unidirectional and bidirectional models in detail for a selected hour and pipeline.

\subsection{Consequences of the non-optimal gas flow direction}
\begin{figure}[h]
\centering
    \begin{subfigure}[b]{0.35\textwidth}
        \centering
        \includegraphics[width=1\linewidth]{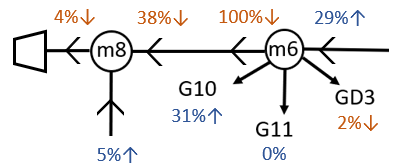}
        \caption{Unidirectional gas flow model}
        \label{fig:LP_zoom} 
    \end{subfigure}
    \begin{subfigure}[b]{0.35\textwidth}
        \centering
        \includegraphics[width=1\linewidth]{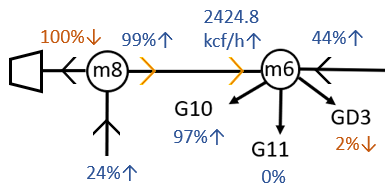}
        \caption{Bidirectional gas flow model}
        \label{fig:MILP_zoom}
    \end{subfigure}
\caption{Close case study of nodes m$8$ and m$6$ showing the percentage changes in gas flows, gas consumption, and production levels in hour $15$ with respect to hour $14$. The percentage change in the production of GFPPs is computed with respect to their installed capacity. The gas inflow into node $6$ from $8$ in plot \eqref{fig:MILP_zoom} is not represented as a percentage, because the inflow at hour $14$ is zero.}
\label{fig:consequence}
 \vspace{-0.2cm}
\end{figure}

To closely examine the reason for the reduced operational cost in the bidirectional gas flow model, we zoom into hour $15$ in the meshed system. Fig. \ref{fig:consequence} focuses on the pipeline from node $m6$ to $m8$, and depicts the percentage change in gas flows and gas load at hour $15$ with respect to the previous hour. In contrast to the pipeline flows, the percentage change in the production of GFPPs G$10$ and G$11$ is computed with respect to their installed capacity. 

In hour 15, the electricity consumption increases by $15.8$\% with respect to the previous hour. The cheapest available power generator to meet the increased demand is GFFP G$10$.
In the unidirectional gas flow model,  the predetermined gas flow direction from node $6$ to $8$ restricts the amount of natural gas available for the operation of G$10$. As a consequence, G$10$ is only able to ramp up $31$\% of its capacity in the unidirectional gas flow model compared to $97$\% in the bidirectional case. Hence, more expensive generators have to be dispatched in the unidirectional gas flow model, resulting in a higher total operational cost.
Note that neither ramping nor unit commitment constraints have been included in the model, which might impose additional technical constraints on the operational flexibility of GFPPs. 

\subsection{Linepack flexibility}
\begin{figure}[h]
    \centering
    \begin{subfigure}[b]{0.3\textwidth}
        \centering
        \includegraphics[width=\textwidth]{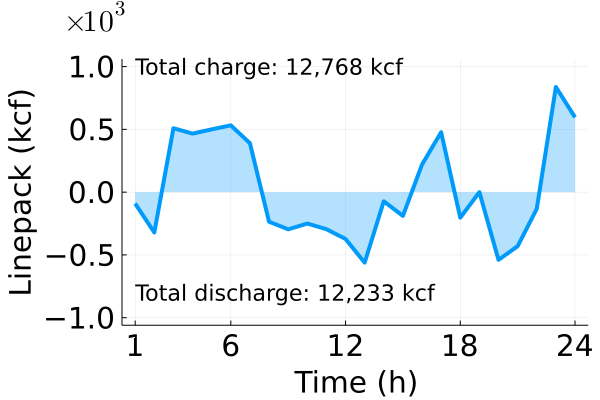}
        \subcaption{Unidirectional, ($m6-m8$)}
        \label{fig:LP pipeline (m6,m8)}
    \end{subfigure}
    \begin{subfigure}[b]{0.3\textwidth}  
        \centering 
        \includegraphics[width=\textwidth]{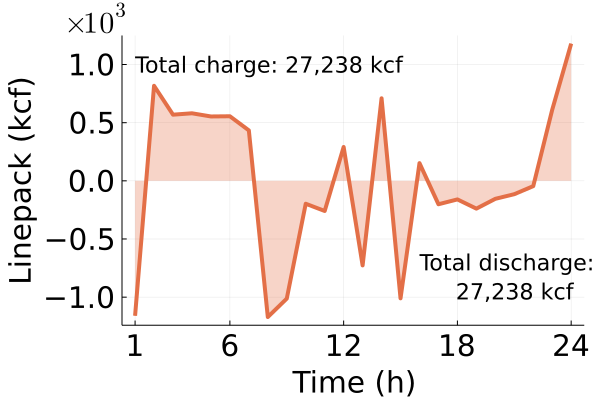}
        \caption{Bidirectional, ($m6-m8$)}  
        \label{fig:MILP pipeline (m6,m8)}
    \end{subfigure}
    \vskip\baselineskip
    \begin{subfigure}[b]{0.3\textwidth}   
        \centering 
        \includegraphics[width=\textwidth]{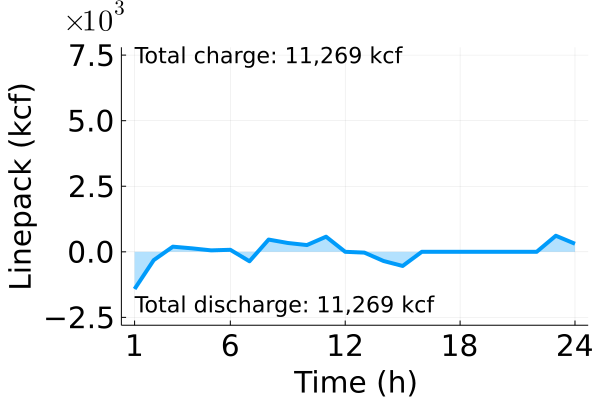}
        \caption{Unidirectional, ($m10-m11$)}    
        \label{fig:LP pipeline (m10,m11)}
    \end{subfigure}
    \begin{subfigure}[b]{0.3\textwidth}   
        \centering 
        \includegraphics[width=\textwidth]{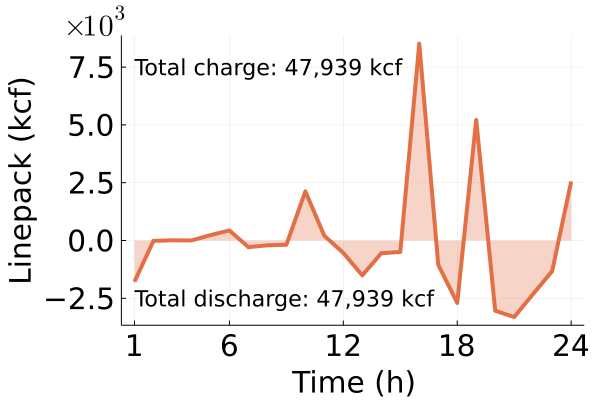}
        \caption{Bidirectional, ($m10-m11$)} 
        \label{fig:MILP pipeline (m10,m11)}
    \end{subfigure}
    \caption{Linepack mass in pipelines ($m6-m8$) and ($m10-m11$) for unidirectional and bidirectional gas flow models.}
    \label{fig:linepack}
\end{figure}

Linepack provides additional operational flexibility to the power system based on the temporal separation of gas inflows and outflows of pipelines due to slow flow transients. Fig.~\ref{fig:linepack} shows the total charge and discharge in two pipelines ($m6-m8$) and ($m10-m11$). We select these two pipelines, since we observe frequent changes of flow direction in such pipelines. A sample of these changes has already been illustrated and discussed in Fig. \ref{fig:consequence} for the case of meshed network. 
We observe that the magnitude of total charge and discharge in the bidirectional gas flow model is comparatively higher than that in the unidirectional gas flow model. This is not necessarily the case for each individual pipeline, but in overall, the bidirectional model charges and discharges more than the unidirectional model.
This indicates that the optimal determination of gas flow directions by treating them as state variables optimizes the charging and discharging capabilities in the pipelines, increasing the overall flexibility potential of linepack. Apart from decreasing the total operational cost of the integrated system, we hypotesize that the increased utilization of linepack may flatten the nodal gas prices, since cheaper gas sources can be utilized more efficiently.

\subsection{Approximation error}
\begin{figure}[h]
\centering
\includegraphics[width=0.6\textwidth,height=6cm]{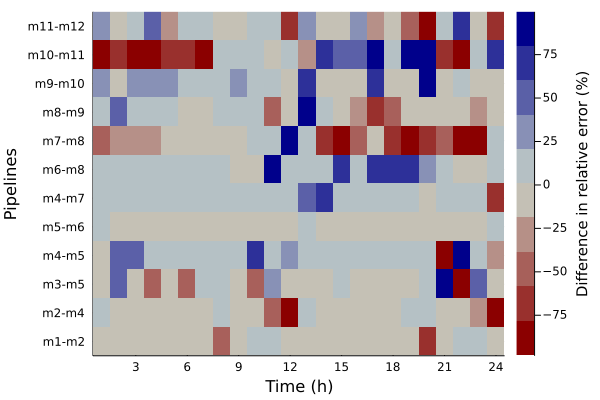}
\caption{Difference in relative errors between left-hand side and right-hand side of Weymouth equation. The blue (red) colors indicate a higher (lower) approximation error in the unidirectional compared to the bidirectional gas flow model.}
\label{fig:approx_error}
\end{figure}

In order to validate our results, we look at the difference in the approximation error for bidirectional and unidirectional gas flow models. The approximation error $\Delta_{m,u,t}$ is calculated as the normalized absolute difference of the optimal natural gas flows of the convexified problem and the flows obtained when using the original Weymouth equation \eqref{eq:MILP_Weymouth} based on pressures. The superscript * indicates optimal values obtained from the models.
\begin{gather}
\Delta_{m,u,t} = \frac{|q_{m,u,t}^{*2}-K_{m,u}^2(pr_{m,t}^{*2}-pr_{u,t}^{*2})|}{K_{m,u}^2(pr_{m,t}^{*2}-pr_{u,t}^{*2})},\hspace{2mm}\forall (m,u) \in \mathcal{Z},t. \label{eq:Delta} 
\end{gather}
Fig. \ref{fig:approx_error} shows the difference in the normalized approximation error in percentage. The blue scale indicates to what extent the bidirectional gas flow model accurately approximates the Weymouth equation in comparison to the unidirectional model, and vice versa for the red scale. We  observe large differences for individual pipelines and time periods. To quantify the overall difference between both models, we use the normalized root mean square error as
\begin{gather}
\Xi = \left[\frac{1}{|\mathcal{T}|\cdot |\mathcal{Z}|} \sum_{t \in \mathcal{T}} \sum_{(m,u) \in \mathcal{Z}} \Delta_{m,u,t}^2 \right]^{\frac{1}{2}}. \label{eq:Xi} 
\end{gather}
where $|\mathcal{T}|$ and $|\mathcal{Z}|$ are the number of time periods and pipelines, respectively. In our case study, $|\mathcal{T}|=24$ and $|\mathcal{Z}|=12$. This error is $0.640$ and $0.636$ for the bidirectional and unidirectional models, respectively, implying that the overall error in both models induced by the approximation of the Weymouth equation is roughly the same. We also notice that the approximation accuracy is highly dependent on the number and choice of fixed pressure points used for the Taylor-series expansion of the Weymouth equation.
 \\
Improving the tightness of the approximation is critical to get meaningful flow levels and validity of the model. For this paper, we solely present the comparison between the errors computed. We find out that the bidirectional model does neither improve nor worsen the overall approximation of the Weymouth equation, compared to the unidirectional model.

%% file: 5-Conclusion.tex
\section{Conclusion and future work}
\label{sec: Section 5}
We provided a detailed analysis of the increased flexibility that the natural gas system provides to the power system when gas flow directions are considered as state variables in the co-optimization problem. This increased flexibility was quantified in terms of the reduced operational cost of the integrated system, by comparing the optimal cost achieved in models with unidirectional and bidirectional gas flows. Using a stylized case study, we found out that modeling gas flows as state variables reduces the operational cost not only in meshed but also in radial networks. The latter particularly happens in a radial network where natural gas sources are located far away from each other in the network.

Given a set of predefined fixed pressure points, we convexified gas flow equations using a Taylor-series expansion. It is worth mentioning that this approximation is usually not tight at optimum, and therefore it may result in schedules which are not necessarily feasible in terms of gas flow dynamics. Consequently, the approximation may overestimate the flexibility that does not exist in the real-time operation. To avoid it, we improved the currently prevalent approximation by tightening the resulting feasible region, ensuring correct flow directions which are consistent with the original non-convex Weymouth equation. Furthermore, we quantified the error of approximation by comparing the convexified and the original Weymouth equations. We found out that the overall errors including all pipelines and time periods in the unidirectional and bidirectional models are in the same order of magnitude. There could be, however, large differences when it comes to individual pipelines and time periods.

As potential directions for the future work, the impact of the approximation error and its relation to the flexibility potential should be further analyzed. One may also conduct a fairness analysis, exploring how the approximation error impacts every individual in terms of the associated operational cost (or profit in a market context). In parallel, an improvement in convexification methods is necessary. Possible approaches and directions have been extensively explored in other research fields, e.g., for the alternating current (AC) power flow problem \cite{Venzke2020InexactFeasibility,Kyri}. The outer approximation method based on the Taylor-series expansion can be further improved by investigating how to more efficiently select the  fixed pressure points as the input data. Furthermore, alternative approaches, e.g., based on machine learning \cite{Bertsimas2021OnlineMilliseconds}, can be explored to reduce the computational time of the resulting MILP in the bidirectional model, by predicting the inactive constraints and/or the value of  binary variables associated with the flow directions at the optimal point. 

This paper exploited a simplified model for the representation of compressors, while discarding other critical gas network components such as valves. This simplification certainly affects the flexibility potential to be unlocked by introducing the gas flow directions as state variables. It is of interest to leverage more detailed operational models for compressors and valves \cite{misha}. The future work should also model potential sources of uncertainty, e.g., renewable power supply, and explores how uncertainty is being propagated from power to natural gas system, depending on the topology of both networks and their interconnection. Recall that the proposed model is a co-optimization, which provides an ideal benchmark. However, power and gas systems are not necessarily being operated by the same entities in the real world. Therefore, the future work should explore how the maximum flexibility potential can be efficiently unlocked in practice \cite{pascal}.  Finally, it is of importance to focus on the transition towards the large deployment of distributed energy resources, in particular electrolyzers and fuel cells, and their associated impacts on the operation of integrated system. 

%% file: 6-Appendix.tex
\section*{Appendices}
\subsection{Objective function \eqref{compact1}}
\label{Appendix:obj}
The objective function minimizes the total operational cost of the integrated system, including the cost of non-GFPPs as well as gas suppliers:
\begin{flalign}
    \label{objective}
    &\underset{\mathrm{\Theta}}{\operatorname{Minimize}} \quad \sum_{t \in \mathcal{T}} \Big( \sum_{i \in \mathcal{C}} C_i^{\rm{E}} p_{i,t} + \sum_{k \in \mathcal{K}} C_k^{\rm{G}} g_{k,t} \Big),
\end{flalign}
where $p_{i, t}$ and $g_{k, t}$ denote the schedules of power generator $i$ and gas supplier $k$ in time period $t$, with associated production and supply costs $C_i^{\rm{E}}$ and $C_k^{\rm{G}}$, respectively. Furthermore, $\mathcal{C}$ denotes the set of non-GFPPs. Recall that $\mathcal{K}$ is the set of gas suppliers. The set of primal variables is $\Theta_{\mathrm{Unidirectional}}$ in the unidirectional gas flow model and $\Theta_{\mathrm{Bidirectional}}$ in the bidirectional case, as already defined in Table \ref{Table:LPvsMILP}.

\subsection{Power system constraints \eqref{compact2} and \eqref{compact3}}
\label{Appendix:power}
Let $(n,r)\in \mathcal{L}$ denote the set of power system lines. Using linearized lossless power flow equations, we enforce power system constraints by
\begin{subequations}\label{power constraints}
\begin{align}
    &0 \leq p_{i,t} \leq P^{\text{max}}_i , \hspace{2mm} \forall i,t \label{LP:p1}\\
    &0 \leq w_{j,t} \leq W_{j,t}, \hspace{2mm} \forall j,t \label{LP:p2}\\
    &f_{n,r,t} = B_{n,r} (\theta_{n,t}-\theta_{r,t}), \hspace{2mm} \forall(n,r) \in \mathcal{L},t \label{LP:p3}\\
    &-F^{\text{max}}_{n,r} \leq f_{n,r,t} \leq F^{\text{max}}_{n,r}, \hspace{2mm} \forall(n,r) \in \mathcal{L},t \label{LP:p4}\\
    &-\pi \leq \theta_{n,t} \leq \pi, \hspace{2mm}\forall n,t \label{LP:p5}\\
    &\theta_{n,t} = 0, \hspace{2mm}\forall n:ref,t \label{LP:p6} \\
    &\sum_{i \in \mathcal{A}^{\rm{I}}_n} p_{i,t} + \sum_{j \in \mathcal{A}^{\rm{J}}_n} w_{j,t} - \sum_{(n,r) \in \mathcal{L}} f_{n,r,t} = \sum_{l \in \mathcal{A}^{\rm{DE}}_n} D^{\mathrm{E}}_{l,t}, \hspace{2mm} \forall n,t. \label{LP:p7} 
\end{align}
\end{subequations}

Constraints \eqref{LP:p1} limit the power production schedule $p_{i,t}$ of generator $i$ in time period $t$ to its installed capacity $P^{\text{max}}_i$. The production schedule of wind farm $j$ in time period $t$, denoted by $w_{j,t}$, is restricted in \eqref{LP:p2} by its single-point deterministic forecast $W_{j,t}$. 
Constraints \eqref{LP:p3} compute the power flow $f_{n,r,t}$ across the line connecting nodes $n$ and $r$ in time period $t$ as the product of the line susceptance $B_{n,r}$ and the difference of the nodal voltage angles $\theta_{n,t}$ and $\theta_{r,t}$. Constraints \eqref{LP:p4} enforce transmission  capacity constraints, where $F^{\text{max}}_{n,r}$ is the capacity of line connecting nodes $n$ and $r$. Constraints \eqref{LP:p5} limit nodal voltage angles, whereas \eqref{LP:p6} fixes the voltage angle at the reference node to zero. The power balance in every node $n$ is enforced by \eqref{LP:p7}, ensuring that the power demand $D^{\mathrm{E}}_{l,t}$ of all loads $l$ located at node $n$ is fully met in each time period $t$. The set of power generators, wind farms, and electricity loads in node $n$ is denoted by $\mathcal{A}^{\rm{I}}_n$, $\mathcal{A}^{\rm{J}}_n$, and $\mathcal{A}^{\rm{DE}}_n$, respectively.

\subsection{Linear approximation of Weymouth equation}
\label{Appendix:Weymouth}
To convexify the non-convex Weymouth equation \eqref{eq:Weymouth} and \eqref{eq:MILP_Weymouth}, it has to be first relaxed and divided into two parts to represent each gas flow direction:
\begin{subequations}\label{eq:weymouth_split}
\begin{align}
    &q_{m,u,t} \leq K_{m,u} \sqrt{pr^2_{m,t} - pr^2_{u,t}}, \hspace{2mm}\forall (m,u) \in \mathcal{Z},t \label{eq:q_mu} \\
    &q_{u,m,t} \leq K_{m,u} \sqrt{pr^2_{u,t} - pr^2_{m,t}}, \hspace{2mm}\forall (m,u) \in \mathcal{Z},t. \label{eq:q_um}
\end{align}
\end{subequations}

Afterwards, a Taylor series approximation is performed around a set of given fixed pressure points $(PR_{m,v},PR_{u,v})$ to approximate the non-convex equations \eqref{eq:weymouth_split} \cite{Tomasgard2007OptimizationChain,Ordoudis2019}. The tightness of this approximation can be controlled by the number of pressure points $v$ at the cost of increasing computational complexity.

Since the flow directions are predetermined in the unidirectional gas flow model, we only need to consider the flow \eqref{eq:q_mu}, which is approximated by
\begin{gather}
    q_{m,u,t} \leq \frac{K_{m,u} PR_{m,v}}{\sqrt{PR^2_{m,v}-PR^2_{u,v}}}pr_{m,t} - \frac{K_{m,u} PR_{u,v}}{\sqrt{PR^2_{m,v}-PR^2_{u,v}}}pr_{u,t}, \hspace{2mm}\forall (m,u) \in \mathcal{Z},v,t. \label{eq:LP_weymouth}
\end{gather}

In the bidirectional gas flow model, the flow directions are endogenously determined in the co-optimization model. Thus, we approximate \eqref{eq:q_mu} and \eqref{eq:q_um} by properly accounting for the flow direction as
\begin{subequations} \label{MILP_flow_approx}
\begin{align}
    &q_{m, u, t}^{+} \leq \frac{K_{m,u} PR_{m,v}}{\sqrt{PR^2_{m,v}-PR^2_{u,v}}}pr_{m,t} - \frac{K_{m,u} PR_{u,v}}{\sqrt{PR^2_{m,v}-PR^2_{u,v}}}pr_{u,t}
    \hspace{2mm}+ M\left(1-y_{m, u, t}\right), \hspace{2mm}\forall\{(m, u) \in \mathcal{Z}| m>u\},v, t \label{MILP6}\\ 
    &q_{m, u, t}^{-} \leq \frac{K_{m,u} PR_{u,v}}{\sqrt{PR^2_{u,v}-PR^2_{m,v}}}pr_{u,t} - \frac{K_{m,u} PR_{m,v}}{\sqrt{PR^2_{u,v}-PR^2_{m,v}}}pr_{m,t}
    \hspace{2mm}+M y_{m, u, t},\hspace{2mm} \forall\{(m, u) \in \mathcal{Z}| m<u\},v, t. \label{MILP7}
\end{align}
\end{subequations}

The terms including the binary variables $M(1-y_{m, u, t})$ and $M y_{m, u, t}$ ensure that when either $q_{m, u, t}^{+}$ or $q_{m, u, t}^{-}$ takes a zero value, the corresponding constraint \eqref{MILP6} or \eqref{MILP7} is inactive.